\documentclass[12pt]{article}

\usepackage{amsfonts}
\usepackage{amsmath}
\usepackage{amsthm}
\usepackage{amssymb}
\usepackage{graphicx}
\usepackage{here}
\usepackage[subrefformat=parens]{subcaption}

\numberwithin{equation}{section}
\usepackage{color}
\definecolor{DarkBlue}{rgb}{0,0,0.7}

\definecolor{DarkRed}{rgb}{0.65,0,0}

\definecolor{DarkGreen}{rgb}{0,0.3,0}

\definecolor{purple}{rgb}{0.7,0,0.7}

\usepackage{ulem}

\usepackage{ascmac}
\usepackage{color}

\allowdisplaybreaks

\theoremstyle{plain}
\newtheorem{main}{}[section] 
\newtheorem{thm}[main]{Theorem}
\newtheorem{lem}[main]{Lemma}

\newtheorem{prop}[main]{Proposition}
\theoremstyle{definition}
\newtheorem{defn}[main]{Definition}
\theoremstyle{remark}
\newtheorem{rem}[main]{Remark}

\begin{document}
% \maketitle

\rightline{\baselineskip16pt\rm\vbox to20pt{
			{%\hbox{APCTP-Pre2010-xxx}
				% \hbox{OCU-PHYS-XXX}
				\hbox{AP-GR-208}
			}
			\vss}}%

	\begin{center}
		{\LARGE Three-Dimensional Almost Contact Metric Manifolds Revisited via the Newman-Penrose Formalism}\\
		\bigskip\bigskip
		{\large
			Satsuki Matsuno\footnote{Satsuki\_Matsuno@omu.ac.jp,\ CRCID:0000-0003-0841-8033}
			% Fumihiro Ueno\footnote{ueno-fumihiro-dx@alumni.osaka-u.ac.jp
		}\\
		\bigskip
		{\it Department of Physics, Graduate School of Science,
			Osaka Metropolitan University\\
			3-3-138 Sugimoto, Sumiyoshi, Osaka 558-8585, Japan}
	\end{center}

\begin{abstract}
This paper applies the Newman-Penrose formalism—a technique primarily used in General Relativity—to the analysis of three-dimensional almost contact metric (ACM) manifolds.
We reformulate and discuss several known notions and properties within the Newman-Penrose framework, demonstrating the applicability of the method in this geometric context.
Furthermore, as an application of the formalism, we address the classification of three-dimensional compact normal ACM manifolds, or equivalently trans-Sasakian manifolds, that admit an $\eta$-Einstein metric.
\end{abstract}

\tableofcontents

\section{Introduction}
The Newman-Penrose formalism has been a powerful tool for the analysis of four-dimensional spacetime \cite{Geroch:1973am}.
It has been adapted to three-dimensional spacetime, and more recently, applied to three-dimensional Riemannian manifolds \cite{Aazami}.
In this paper, we apply the Newman-Penrose formalism to three-dimensional almost contact metric (ACM) manifolds.
The original Newman-Penrose formalism in four-dimensional spacetime has been particularly effective in analyzing the geometry of shear-free geodesic null congruences.
Since three-dimensional normal ACM manifolds can be characterized by the property that the characteristic vector field generates a shear-free geodesic congruence, this observation establishes a clear analogy regarding the context in which the Newman-Penrose formalism can demonstrate its power.

Our primary interest lies in analyzing the geometry of three-dimensional ACM manifolds using the Newman-Penrose formalism.
The definitions of various classes of ACM manifolds proposed to date—such as contact metric manifolds and trans-Sasakian manifolds—do not always make the underlying geometric situation immediately clear.
By employing the Newman-Penrose formalism, we can describe these structures and their properties in terms of basic geometrical quantities—namely, twist, expansion, and shear—which are fundamental properties of geodesic congruences.
In this framework, we reformulate several definitions and known properties of three-dimensional ACM manifolds using the language of the Newman-Penrose formalism.

Furthermore, to demonstrate the effectiveness of the Newman-Penrose formalism in the analysis of three-dimensional ACM manifolds, we consider $\eta$-Einstein metrics, in particular the classification of three-dimensional compact normal ACM manifolds that admit an $\eta$-Einstein metric, where in dimension three a normal ACM structure is also referred to as a trans-Sasakian structure.
The case of three-dimensional contact metric structures has already been clarified by Blair et al. \cite{Blair1990}.
Meanwhile, the necessary and sufficient conditions for a three-dimensional normal ACM manifold to be $\eta$-Einstein were reported by Olszak \cite{Olszak}.
Subsequent research has further elucidated several equivalent conditions, geometric properties, and their relation to Ricci solitons \cite{Marrero,De2003,De2008,Al-Solamy,Ganguly}.
However, to the best of our knowledge, we have been unable to locate any prior literature explicitly reporting classification results for three-dimensional compact normal ACM structures subject solely to the $\eta$-Einstein condition.
By introducing the Newman-Penrose formalism, we find that the $\eta$-Einstein condition possesses a non-linear subelliptic structure.
When the manifold is compact, we can establish certain constancy results using a priori estimates.
This provides a classification result of $\eta$-Einstein manifolds for compact three-dimensional normal ACM manifolds.

\textbf{Classification Result}:
Let \( (M,\varphi,\xi,\eta,g) \) be a compact three-dimensional normal almost contact metric manifold, or equivalently a trans-Sasakian manifold. If \( (M,\varphi,\xi,\eta,g) \) is $\eta$-Einstein, then the structure is either
(i) a Sasakian structure up to homothety, or (ii) a $\beta$-Kenmotsu structure whose function $\beta$ varies only along $\xi$.

The organization of this paper is as follows.
In Section \ref{sec:preliminary}, we summarize the Newman-Penrose formalism on three-dimensional Riemannian manifolds.
In Section \ref{sec:reformulation}, we reformulate the three-dimensional ACM structure using the Newman-Penrose formalism.
Several traditional definitions in the context of ACM geometry are expressed in the language of the Newman-Penrose formalism.
In Section \ref{sec:eta-einstein problem}, we address the $\eta$-Einstein condition on three-dimensional contact metric structures and normal ACM structures.
Although the results for three-dimensional contact metric manifolds are known, we provide an alternative proof using the Newman-Penrose method. Subsequently, we prove the main result using the Newman-Penrose formalism.

\section{Preliminary}\label{sec:preliminary}
In this section, we summarize the Newman-Penrose formalism on three-dimensional Riemannian manifolds.
The material in this section is based on \cite{Geroch:1973am,Aazami}.
Let $(M,g)$ be an orientable three-dimensional Riemannian manifold.
We can choose a nowhere vanishing unit vector field $\xi$ and a suitable orthonormal frame $\{\xi,e_2,e_3\}$ at least locally.
Setting $\partial=\frac{1}{\sqrt{2}}(e_2-ie_3)$, we obtain a local complex frame $\{\xi,\partial,\bar{\partial}\}$.
In this setting, the metric satisfies $g(\xi,\xi)=g(\partial,\bar{\partial})=1$ and $g(\xi,\partial)=g(\partial,\partial)=g(\bar{\partial},\bar{\partial})=0$.
Letting $\{\eta,\theta^2,\theta^3\}$ be the dual frame of $\{\xi,e_2,e_3\}$ and setting $\mu=\frac{1}{\sqrt{2}}(\theta^2+i\theta^3)$, we see that the complex dual frame corresponding to the complex frame $\{\xi,\partial,\bar{\partial}\}$ is $\{\eta,\mu,\bar{\mu}\}$ by noting that $\mu(\partial)=1$.
The \textbf{spin coefficients} with respect to the frame $\{\xi,\partial,\bar{\partial}\}$ are defined as follows.
\begin{align}
\kappa=-g(\nabla_\xi\xi,\partial),\ \rho=g(\nabla_{\bar{\partial}}\xi,\partial),\ \sigma=-g(\nabla_\partial\xi,\partial),\\
\epsilon=g(\nabla_\xi\partial,\bar{\partial}),\ \beta=g(\nabla_\partial\partial,\bar{\partial}).
\end{align}
The above definition of $\rho$ differs in sign from the one in \cite{Aazami}. This choice is convenient in the context of contact metric geometry.
These are equivalently defined as the coefficients of the following covariant derivatives or Lie brackets.
\begin{align}
&\nabla_\xi \xi = -\bar{\kappa}\partial - \kappa\bar{\partial},\ 
\nabla_\partial \xi = \bar{\rho}\partial - \sigma\bar{\partial},\ 
\nabla_\xi \partial = \kappa\xi + \epsilon\partial,\\
&\nabla_\partial\partial = \sigma\xi + \beta\partial,\ 
\nabla_{\partial} \bar{\partial} = -\bar{\rho}\xi - \beta\bar{\partial},\ \\
&[\xi,\partial] = \kappa\xi + (\epsilon - \bar{\rho})\partial + \sigma\bar{\partial},\ 
[\partial,\bar{\partial}] = (\rho - \bar{\rho})\xi + \bar{\beta}\partial - \beta\bar{\partial}.
\end{align}

If $\kappa = 0$, the integral curves of $\xi$ are geodesics.
We write $\rho = \Theta + i\omega$.
The real part $\Theta$ is called the expansion of $\xi$, and the imaginary part $\omega$ is called the twist of $\xi$.
The coefficient $\sigma$ is called the complex shear of $\xi$.
The quantities $\Theta$, $\omega$, and $\sigma$ admit the following geometric interpretations.
Specifically, considering the congruence generated by $\xi$, $\Theta$ represents the rate of change of the area of a cross-section orthogonal to $\xi$, $\omega$ represents the angular velocity of the cross-section, and $\sigma$ represents the shear deformation of the cross-section.
The coefficient $\epsilon$ is purely imaginary, and one can locally choose a frame such that $\epsilon$ vanishes.
The coefficient $\beta$ is related to the sectional curvature of the plane orthogonal to $\xi$.
On a homogeneous manifold, there often exists a frame such that $\beta=0$.

The choice of a complex vector field $\partial$ is not unique, or $\partial$ may not be chosen globally and hence may have to be replaced somewhere.
Let $U\subset M$ be a open neighborhood and $e^{i\theta}: U \to S^1$ be a smooth map.
Then a gauge transformation is given by $\partial\to\partial'=e^{i\theta}\partial$.
Under this gauge transformation, the spin coefficients transform as follows:
\begin{align}
&\kappa'=e^{i\theta}\kappa,\ \rho'=\rho,\ \sigma'=e^{2i\theta}\sigma,\label{eq:spin coeff trans}\\
&\epsilon'=\epsilon+i\xi(\theta),\ \beta'=e^{i\theta}(\beta+i\partial\theta).
\end{align}
As can be seen from these transformations, the value of $\rho$, namely the expansion and twist, is well-defined.
However, for $\kappa$ and $\sigma$, only statements concerning their magnitudes, $|\kappa|$ and $|\sigma|$, are well-defined.
A quantity $q$ is said to be a \textbf{spin-weighted quantity} of \textbf{spin weight} $s$ if it transforms as $q'=e^{is\theta}q$ under a gauge transformation.
The spin weights of $\kappa$, $\rho$, and $\sigma$ are $1$, $0$, and $2$, respectively.
Particular attention must be paid to the differentiation of spin-weighted quantities; the derivatives $\xi(q)$, $\partial(q)$, and $\bar{\partial}(q)$ are not, in general, spin-weighted quantities themselves.
However, the following operators map spin-weighted quantities to spin-weighted quantities \cite{Geroch:1973am}:
\begin{align}
P(q)&=\xi(q)-s\epsilon q,\\
\eth(q)&=\partial(q)-sq\beta,\\
\bar{\eth}(q)&=\bar{\partial}(q)+sq\bar{\beta}.
\end{align}
where $s$ is the spin weight of $q$.
The operator $P$ is pronounced ``thorn'' and $\eth$ is pronounced ``eth''.
A straightforward calculation verifies that $P(q)$, $\eth(q)$, and $\bar{\eth}(q)$ are spin-weighted quantities of spin weights $s$, $s+1$, and $s-1$, respectively, as detailed in Appendix \ref{App: P eth spin covariance}.
Noting that complex conjugation reverses the sign of the spin weight, we obtain the relation $\overline{\bar{\eth}(q)}=\eth(\bar{q})$.

The spin coefficients are not mutually independent.
The \textbf{generalized Sachs equations} are given as follows \cite{Aazami}.
\begin{align}
&-\xi(\rho) - \bar{\partial}(\kappa)
    = |\kappa|^{2} + |\sigma|^{2} + (\rho)^{2}
      + \kappa \bar{\beta} + \tfrac12 \mathrm{Ric}(\xi,\xi),\label{eq:1st bianchi 1}\\
&\xi(\sigma) - \partial(\kappa)
    = (\kappa)^{2} + 2\sigma\epsilon
      - \sigma(\rho + \bar{\rho})
      - \kappa\beta + \mathrm{Ric}(\partial,\partial),\label{eq:1st bianchi 2}\\
&-\partial(\rho) - \bar{\partial}(\sigma)
    = 2\sigma\bar{\beta} + (\rho - \bar{\rho})\kappa
      + \mathrm{Ric}(\xi,\partial),\label{eq:1st bianchi 3}\\
&\xi(\beta) - \partial(\epsilon)
    = \sigma(\bar{\kappa} - \bar{\beta})
      + \kappa(\epsilon + \bar{\rho})
      + \beta(\epsilon - \bar{\rho})
      - \mathrm{Ric}(\xi,\partial),\label{eq:1st bianchi 4}\\
&\partial(\bar{\beta}) + \bar{\partial}(\beta)
    = |\sigma|^{2} - |\rho|^{2} - 2|\beta|^{2}
      - (\rho - \bar{\rho})\epsilon
      - \mathrm{Ric}(\partial,\bar{\partial})
      + \tfrac12 \mathrm{Ric}(\xi,\xi).\label{eq:1st bianchi 5}
\end{align}
The first three equations of the generalized Sachs equations can be rewritten using the operators $P$ and $\eth$ in such a way that both sides are spin-weighted quantities.
Thus, these equations are well-defined independently of the choice of a complex frame.
Similarly, the second Bianchi identity yields the following equations.
\begin{align}
&\xi(\text{Ric}(\xi, \partial)) - \frac{1}{2}\partial(\text{Ric}(\xi, \xi)) + \overline{\partial}(\text{Ric}(\partial, \partial)) \\
&= \kappa\text{Ric}(\xi, \xi) + (\epsilon - 2\rho - \bar{\rho})\text{Ric}(\xi, \partial) \\
&+ \sigma\text{Ric}(\xi, \overline{\partial}) - (\overline{\kappa} + 2\overline{\beta})\text{Ric}(\partial, \partial) - \kappa\text{Ric}(\partial, \overline{\partial}),\label{eq:2nd bianchi 1}\\
&\partial(\text{Ric}(\xi, \overline{\partial})) + \overline{\partial}(\text{Ric}(\xi, \partial)) - \xi(\text{Ric}(\partial, \overline{\partial})) + \frac{1}{2}\xi(\text{Ric}(\xi, \xi)) \\
&=-(\rho + \overline{\rho}) (\text{Ric}(\xi, \xi) - \text{Ric}(\partial, \bar{\partial})) \\
&- \overline{\sigma}\text{Ric}(\partial, \partial) - \sigma\text{Ric}(\overline{\partial}, \overline{\partial}) - (2\overline{\kappa} + \overline{\beta})\text{Ric}(\xi, \partial) - (2\kappa + \beta)\text{Ric}(\xi, \overline{\partial}).\label{eq:2nd bianchi 2}
\end{align}

\section{Characterizations of three-dimensional ACM structures via the Newman--Penrose formalism}\label{sec:reformulation}
In this section we reformulate the definitions of several well-known notions on three-dimensional ACM manifolds, such as normality and contact metric structures, in terms of the spin coefficients.
We also reformulate several well-known classes of manifolds including trans-Sasakian manifolds and $(k,\mu,\nu)$-manifolds, and we discuss some of their known properties.

Let $(M,\varphi,\xi,\eta,g)$ be a three-dimensional ACM manifold.
Here $\xi$ is the characteristic vector field, $\eta$ is the metric dual of $\xi$, and $\varphi\colon\Gamma(TM)\to\Gamma(TM)$ is an endomorphism field satisfying
\begin{align}
&\varphi^2=-I+\eta\otimes\xi,\ \eta(\xi)=1,\ \varphi(\xi)=0,\\   
&g(\varphi(X),\varphi(Y))=g(X,Y)-\eta(X)\eta(Y),\ X,Y\in\Gamma(TM).
\end{align}
Let $\{U_\alpha\}_{\alpha\in\Lambda}$ be an open covering on $M$, and let $\{\partial^{(\alpha)},\bar{\partial}^{(\alpha)}\}$ be a complex frame orthogonal to $\xi$ on $U_\alpha$.
We can always choose $\partial^{(\alpha)}$ so that $\varphi(\partial^{(\alpha)})=i\partial^{(\alpha)}$ holds.
Conversely, we can recover all the information about the ACM structure $(\varphi,\xi,\eta,g)$ from the frame $\{\xi,\partial^{(\alpha)},\bar{\partial}^{(\alpha)}\}$ on $U_\alpha\ (\alpha\in\Lambda)$.
Suppose that a frame $\{\xi,\partial^{(\alpha)},\bar{\partial}^{(\alpha)}\}$ is given on an orientable three dimensional Riemannian manifold $(M,g)$, and $\xi$ is a global unit vector field.
Since the metric $g|_D$ is defined on the distribution $D$ orthogonal to $\xi$, we can uniquely choose a complex structure $J$ on $D$ such that $g|_D(J\cdot,J\cdot)=g|_D(\cdot,\cdot)$ and $J(\partial^{(\alpha)})=i\partial^{(\alpha)}$ on $U_\alpha$.
Defining an endomorphism field $\varphi$ by $\varphi|_D:=J$ and $\varphi(\xi)=0$, we obtain a unique ACM structure.
Thus we may discuss ACM structures on $M$ in terms of the frame $\{\xi,\partial^{(\alpha)},\bar{\partial}^{(\alpha)}\}$ and an open covering $\{U_\alpha\}_{\alpha\in\Lambda}$ on orientable three-dimensional Riemannian manifolds.
We call this complex frame $\{\xi,\partial^{(\alpha)},\bar{\partial}^{(\alpha)}\}$ on $U_\alpha$ a \textbf{local adapted complex frame}.
The spin coefficients and all other quantities and definitions are the same as in Section \ref{sec:preliminary}.

% Let $(M,\varphi,\xi,\eta,g)$ be a three-dimensional ACM manifold.
% Here $\xi$ is the characteristic vector field, $\eta$ is the metric dual of $\xi$, and $\varphi\colon\Gamma(TM)\to\Gamma(TM)$ is an endomorphism field satisfying
% \begin{align}
% &\varphi^2=-I+\eta\otimes\xi,\ \eta(\xi)=1,\ \varphi(\xi)=0,\\   
% &g(\varphi(X),\varphi(Y))=g(X,Y)-\eta(X)\eta(Y),\ X,Y\in\Gamma(TM).
% \end{align}
% Let $\{\partial,\bar{\partial}\}$ be a complex frame orthogonal to $\xi$.
% We can always choose $\partial$ so that $\varphi(\partial)=i\partial$ holds.
% Conversely, we can recover all the information about the ACM structure $(\varphi,\xi,\eta,g)$ from the frame $\{\xi,\partial,\bar{\partial}\}$.
% Suppose that a frame $\{\xi,\partial,\bar{\partial}\}$ is given on an orientable three dimensional Riemannian manifold $(M,g)$.
% Since the metric $g|_D$ is defined on the distribution $D$ orthogonal to $\xi$, we can uniquely choose a complex structure $J$ on $D$ such that $g|_D(J\cdot,J\cdot)=g|_D(\cdot,\cdot)$ and $J(\partial)=i\partial$.
% Defining an endomorphism field $\varphi$ by $\varphi|_D:=J$ and $\varphi(\xi)=0$, we obtain a unique ACM structure.
% Thus we may discuss ACM structures on $M$ in terms of the frame$\{\xi,\partial,\bar{\partial}\}$ on orientable three-dimensional Riemannian manifolds.
% The spin coefficients and all other quantities and definitions are the same as in Section \ref{sec:preliminary}.

On an odd-dimensional manifold, an ACM structure $(\varphi,\xi,\eta,g)$ is said to be normal if the Nijenhuis tensor of the almost complex structure defined on $M\times\mathbb{R}$ by $\varphi,\xi,\eta$ vanishes.
A necessary and sufficient condition for an ACM structure to be normal is given by
\begin{align}
[\varphi,\varphi]+2\xi\otimes d\eta=0,
\end{align}
where $[\varphi,\varphi]$ denotes the Nijenhuis torsion tensor of $\varphi$.
In dimension three, this condition simplifies significantly, and it has been shown that
\begin{align}
\nabla_{\varphi(X)}\xi=\varphi\nabla_X\xi,
\end{align}
is a necessary and sufficient condition for the structure to be normal \cite{Olszak}.

\begin{screen}
\begin{prop}
Let $(M, \varphi, \xi, \eta, g)$ be a three-dimensional almost contact metric manifold.
The almost contact metric structure is normal if and only if the characteristic vector field $\xi$ generates a shear-free geodesic congruence.
\end{prop}
\end{screen}
\begin{proof}
Let $U$ be a neighborhood and $\{\xi,\partial,\bar{\partial}\}$ be a local adapted complex frame on $U$.
Since the necessary and sufficient condition for a three-dimensional ACM structure to be normal is $\nabla_{\varphi(X)}\xi=\varphi\nabla_X\xi$ and since $\varphi(\partial)=i\partial$, the normality condition is equivalent to the following:
\begin{align}
&\nabla_\xi\xi=0,\\
&i\nabla_{\partial}\xi=\varphi\nabla_{\partial}\xi,
\end{align}
Thus, from the first condition we obtain $\kappa=0$, and from the second condition we obtain
\begin{align}
i\bar{\rho}\partial - i\sigma\bar{\partial}=i\bar{\rho}\partial +i \sigma\bar{\partial},\\
\sigma=0.
\end{align}
By \eqref{eq:spin coeff trans}, the condition $\kappa=\sigma=0$ is independent of the choice of the complex frame $(U,\partial)$, hence the same conclusion holds on all of $M$.
\end{proof}

A three-dimensional ACM structure $(\varphi,\xi,\eta,g)$ is said to be of rank 1 if $d\eta=0$, of rank 2 if $d\eta\ne0$ and $\eta\wedge d\eta=0$, and of rank 3 if $\eta\wedge d\eta\ne0$.
By a direct computation we obtain
\begin{align}
&d\eta=(\kappa\mu+\bar{\kappa}\bar{\mu})\wedge\eta-2i\omega\mu\wedge\bar{\mu},\\
&\eta\wedge d\eta=-2i\omega\eta\wedge\mu\wedge\bar{\mu},
\end{align}
and hence a normal ACM structure cannot be of rank 2.

For an ACM structure $(\varphi, \xi, \eta, g)$, the fundamental two-form is defined by $\Phi(X,Y):=g(X,\varphi(Y))$.
An ACM structure is called a contact metric structure if it satisfies $d\eta=\Phi$.
In dimension three, since $\Phi = -2i\mu\wedge\bar{\mu}$, we readily obtain the following characterization.

\begin{screen}
\begin{prop}
Let $(M, \varphi, \xi, \eta, g)$ be a three-dimensional almost contact metric manifold. 
The almost contact metric structure is a contact metric structure if and only if $\kappa=0$ and $\omega=1$ hold.
Furthermore, on a contact metric manifold, we have $\Theta=0$.
\end{prop}
\end{screen}
The last statement follows from the imaginary part of \eqref{eq:1st bianchi 1}.

From the above proposition, we obtain the following physical interpretation of a contact metric structure: the characteristic vector field $\xi$ generates a geodesic congruence with twist equal to one and expansion-free.
It should be noted that in this case no condition is imposed on the complex shear $\sigma$.

Trans-Sasakian structures provide many examples and have been studied extensively \cite{Marrero,De2003,De2008,Al-Solamy,Ganguly}.
An ACM structure is said to be trans-Sasakian of type $(\alpha,\beta_s)$ if
\begin{align}
(\nabla_X\varphi)Y = \alpha(g(X,Y)\xi - \eta(Y)X) + \beta_s(g(\varphi X, Y)\xi - \eta(Y)\varphi X)\ \ \ \alpha,\beta_s\in C^\infty(M).\label{eq:trans-sasaki def}
\end{align}
In much of the literature, the symbol $\beta$ is used in place of $\beta_s$, but since $\beta$ is also used for a spin coefficient, we shall use $\beta_s$ to avoid confusion, where the subscript $s$ denotes ``Sasakian''.
A trans-Sasakian structure with $\beta_s=0$ is called an $\alpha$-Sasakian structure.
A trans-Sasakian structure with $\alpha=0$ is called a $\beta_s$-Kenmotsu structure.

\begin{screen}
\begin{prop}
Let $(M, \varphi, \xi, \eta, g)$ be a three-dimensional almost contact metric manifold. 
The almost contact metric structure is a trans-Sasakian structure of type $(\alpha,\beta_s)$ if and only if $\kappa=\sigma=0$ and $\rho=\beta_s+i\alpha$.
\end{prop}
\end{screen}
\begin{proof}
Let $U$ be a neighborhood and $\{\xi,\partial,\bar{\partial}\}$ be a local adapted complex frame on $U$.
If we take $X=\xi$ and $Y=\partial$ then we have
\begin{align}
(\nabla_\xi\varphi)(\partial)&=\nabla_\xi(\varphi \partial) - \varphi(\nabla_\xi \partial) = \nabla_\xi(i\partial) - \varphi(\kappa\xi + \epsilon\partial)\\
&= i(\kappa\xi + \epsilon\partial) - (\varphi(\kappa\xi) + \epsilon \varphi(\partial))
= i\kappa\xi,
\end{align}
and hence $\kappa=0$.
Similarly, taiking $X=\partial$ and $Y=\xi$, we compare
\[
\nabla_\partial(\varphi \xi) - \varphi(\nabla_\partial \xi) = -i\bar{\rho}\partial - i\sigma\bar{\partial},
\]
wtih
\[
\alpha(g(\partial,\xi)\xi - \eta(\xi)\partial) + \beta_s(g(\varphi \partial, \xi)\xi - \eta(\xi)\varphi \partial)=-\alpha\partial - i\beta_s\partial = -(\alpha+i\beta_s)\partial,
\]
and we obtain $\sigma=0$ and $\rho=\beta_s+i\alpha$.

Conversely, assume that $\kappa=\sigma=0$ and $\rho=\beta_{s}+i\alpha$.
It is sufficient to verify the condition \eqref{eq:trans-sasaki def} for $(X,Y)=(\partial,\partial)$ and $(X,Y)=(\partial,\overline{\partial})$.
In the case of $(X,Y)=(\partial,\partial)$, the left-hand side becomes
\begin{equation}
    (\nabla_{\partial}\varphi)\partial = \nabla_{\partial}(i\partial) - \varphi(\beta\partial) = i\beta\partial - i\beta\partial = 0,
\end{equation}
and the right-hand side clearly vanishes.
In the case of $(X,Y)=(\partial,\overline{\partial})$, the left-hand side becomes
\begin{align}
    (\nabla_{\partial}\varphi)\overline{\partial} &= \nabla_{\partial}(-i\overline{\partial}) - \varphi(-\overline{\rho}\xi - \beta\overline{\partial}) \\
    &= -i(-\overline{\rho}\xi - \beta\overline{\partial}) - i\beta\overline{\partial} \\
    &= i\overline{\rho}\xi = (\alpha + i\beta_{s})\xi,
\end{align}
and the right-hand side is calculated as
\begin{align}
    \alpha \xi + \beta_{s}(g(i\partial, \overline{\partial})\xi) = (\alpha + i\beta_{s})\xi.
\end{align}
By \eqref{eq:spin coeff trans}, the condition $\kappa=\sigma=0$ and the value of $\rho$ are independent of the choice of the complex frame $(U,\partial)$, hence the same conclusion holds on all of $M$.
\end{proof}

As implied by this proposition, a three-dimensional trans-Sasakian structure is equivalent to a three-dimensional normal ACM structure.
From the above proposition, we can reformulate the definitions of well-known special classes of three-dimensional trans-Sasakian structures as follows.
\begin{screen}
\begin{defn}
Let $(M, \varphi, \xi, \eta, g)$ be a three-dimensional trans-Sasakian manifold of type $(\alpha,\beta_s)$.
The trans-Sasakian structure is called an $\alpha$-Sasakian structure if $(\alpha,\beta_s)=(\omega,0)$.
In particular, an $\alpha$-Sasakian structure of twist $\omega=1$ is called a Sasakian structure.
The trans-Sasakian structure is called a $\beta_s$-Kenmotsu structure if $(\alpha,\beta_s)=(0,\Theta)$.
In particular, a $\beta_s$-Kenmotsu structure of expansion $\Theta=1$ is called a Kenmotsu structure.
The trans-Sasakian structure is called a cosymplectic structure if $(\alpha,\beta_s)=(0,0)$.
\end{defn}
\end{screen}

In dimension three, an $\alpha$-Sasakian structure is also called a quasi-Sasakian structure.
Another notion, called a nearly Sasakian structure, also appears in the literature.
However, in dimension three it reduces to a Sasakian structure, so we do not use this notion in the present paper.
For three-dimensional $\alpha$-Sasakian structures it is well known that $\xi$ is a Killing vector field.
This follows immediately from the following lemma.

\begin{screen}
\begin{lem}
Let $(M, \varphi, \xi, \eta, g)$ be a three-dimensional almost contact metric manifold.
The characteristic vector field $\xi$ is a Killing vector field if and only if $\kappa=\sigma=\Theta=0$.
\end{lem}
\end{screen}
\begin{proof}
This can be readily seen from the following calculation.
Let $U$ be a neighborhood and $\{\xi,\partial,\bar{\partial}\}$ be a local adapted complex frame on $U$.
\begin{align}
-2\sigma&=2g(\nabla_{\partial}\xi,\partial)=g(\nabla_{e_2}\xi,e_2)-g(\nabla_{e_3}\xi,e_3)-i(g(\nabla_{e_2}\xi,e_3)+g(\nabla_{e_3}\xi,e_2)),\\
2\Theta&=2\textrm{Re} (g(\nabla_{\bar{\partial}}\xi,\partial))=g(\nabla_{e_2}\xi,e_2)+g(\nabla_{e_3}\xi,e_3),\\
-\sqrt{2}\kappa&=g(\nabla_\xi\xi,e_2)-ig(\nabla_\xi\xi,e_3).
\end{align}
By \eqref{eq:spin coeff trans}, the condition $\kappa=\sigma=0$ is independent of the choice of the complex frame $(U,\partial)$, hence the same conclusion holds on all of $M$.
\end{proof}

Since a trans-Sasakian structure is normal, the characteristic vector field $\xi$ generates a shear-free geodesic congruence. By contrast, in a contact metric structure the characteristic vector field is geodesic but not necessarily shear-free.
When shear is present, the geometric behavior becomes more complicated and its analysis becomes more difficult.
In the traditional context of contact metric geometry, the tensor $h:=\frac{1}{2}\mathcal{L}_\xi\varphi$ has been used to analyze the shear of a contact metric structure.
To see the relation between $h$ and the complex shear $\sigma$, it suffices to observe the following simple computation.
\begin{align}
2h(\partial) &= [\xi, \varphi(\partial)] - \varphi[\xi, \partial] \\
&= [\xi, i\partial] - \varphi[\xi, \partial] \\
 &= i \left( \kappa\xi + (\epsilon - \overline{\rho})\partial + \sigma\overline{\partial} \right) - \left( 0 + i(\epsilon - \overline{\rho})\partial - i\sigma\overline{\partial} \right) \\
&= i\kappa\xi + i(\epsilon - \overline{\rho})\partial + i\sigma\overline{\partial} - i(\epsilon - \overline{\rho})\partial + i\sigma\overline{\partial}\\
 &= i\kappa\xi + 2i\sigma\overline{\partial}, \\
\therefore \quad h(\partial) &= i\sigma\overline{\partial} + \frac{i}{2}\kappa\xi.
\end{align}
It follows that, when $\kappa=0$, we have
\begin{align}
h(\partial) &= i\sigma\overline{\partial},\\
h(\xi)&=0.
\end{align}
Moreover, the matrix representation of $h$ with respect to the frame $\{e_2,e_3\}$ is given by
$$h = \begin{pmatrix}
-\text{Im}(\sigma) & -\text{Re}(\sigma) \\
-\text{Re}(\sigma) & \text{Im}(\sigma)
\end{pmatrix}.$$
The operator $h$ is real and symmetric, and its eigenvalues are $\pm|\sigma|$.
If $\sigma\neq 0$ and we write $\sigma=|\sigma|e^{i\gamma}$, then the eigendirections corresponding to the eigenvalues $+|\sigma|$ and $-|\sigma|$ are obtained by rotating $e_2$ through the angles $-\frac{\pi}{4}-\frac{\gamma}{2}$ and $\frac{\pi}{4}-\frac{\gamma}{2}$, respectively, modulo $\pi$.

A well-known class that includes contact metric structures with shear is the class of $(k,\mu,\nu)$-structures.
A three-dimensional contact metric manifold is called a $(k,\mu,\nu)$-manifold if it satisfies
\begin{align}
R(X,Y)\xi=k(\eta(Y)X-\eta(X)Y)+\mu(\eta(Y)hX-\eta(X)hY)+\nu(\eta(Y)h\phi X-\eta(X)h\phi Y).
\end{align}
In the traditional notation, $\kappa$ is used instead of $k$, but here we use $k$ in order to avoid confusion with the spin coefficients.
In the three-dimensional Newman-Penrose formalism we have $Ric(\xi,\xi)=2g(R(\partial,\xi)\xi,\bar{\partial})$, $Ric(\partial,\xi)=g(R(\bar{\partial},\partial)\xi,\partial)$, and $Ric(\partial,\partial)=g(R(\xi,\partial)\partial,\xi)$, and hence by a straightforward computation one sees that this condition is equivalent to the following definition.

\begin{screen}
\begin{defn}
Let $(M, \varphi, \xi, \eta, g)$ be a three-dimensional contact metric manifold.
The contact metric structure is called a $(k,\mu,\nu)$-structure if the components of the Ricci tensor are given by
\begin{align}
&\textrm{Ric}(\xi, \xi) = 2k, \\
&\textrm{Ric}(\partial, \xi) = 0, \\
&\textrm{Ric}(\partial,\partial) = (i\mu - \nu)\sigma.
\end{align}
\end{defn}
\end{screen}

Both sides of the equation in the above definition have the same spin weight, and hence the definition is well-defined independently of the choice of the local adapted frame.
Since the definition of a $(k,\mu,\nu)$-structure is a condition on the curvature, it can be expressed in terms of the spin coefficients and their first order derivatives.

\begin{screen}
\begin{prop}\label{prop:(k,mu,nu)-mfd condition}
Let $(M, \varphi, \xi, \eta, g)$ be a three-dimensional contact metric manifold.
The contact metric structure is a $(k,\mu,\nu)$-structure if and only if we have the following equations:
\begin{align}
&\Theta=0,\\
&k = 1 - |\sigma|^{2},\\
&P(\sigma)=\sigma(i\mu-\nu),\\
&\bar{\eth}(\sigma)=0.
\end{align}
\end{prop}
\end{screen}
\begin{proof}
Let $U$ be a neighborhood and $\{\xi,\partial,\bar{\partial}\}$ be a local adapted complex frame on $U$.
Substituting $\kappa=0,\ \rho=\Theta+i$ into the generalized Sachs equation
$$
-\xi(\rho) - \bar{\partial}(\kappa)
    = |\kappa|^{2} + |\sigma|^{2} + (\rho)^{2}
      + \kappa \bar{\beta} + \tfrac12 \mathrm{Ric}(\xi,\xi),
$$
and comparing the real and imaginary parts of both sides, we obtain
\begin{align}
&\Theta=0,\\
&k = 1 - |\sigma|^{2}.
\end{align}
Next, from the generalized Sachs equation
$$
-\partial(\rho) - \bar{\partial}(\sigma)
    = 2\sigma\bar{\beta} + (\rho - \bar{\rho})\kappa
      + \mathrm{Ric}(\xi,\partial),
$$
we obtain
\begin{align}
\bar{\eth}\sigma=\bar{\partial}\sigma+2\sigma\bar{\beta}=0.
\end{align}
Finally, from the generalized Sachs equation
$$
\xi(\sigma) - \partial(\kappa)
    = (\kappa)^{2} + 2\sigma\epsilon
      - \sigma(\rho + \bar{\rho})
      - \kappa\beta + \mathrm{Ric}(\partial,\partial),
$$
we obtain
\begin{align}
P(\sigma)=\sigma(i\mu-\nu).
\end{align}
Conversely, in a contact metric structure, if these conditions on the spin coefficients hold, then from the above three 1st Bianchi identities one can verify that the components of the Ricci tensor appearing in the definition of a $(k,\mu,\nu)$-structure are recovered.

The claimed four conditions are independent of the choice of the complex frame $(U,\partial)$, hence the same conclusion holds on all of $M$.
\end{proof}

In a $(k,\mu,\nu)$-structure we have $k=1$ if and only if $\sigma=0$, and hence in this case the structure is Sasakian.
From this observation one may regard a $(k,\mu,\nu)$-structure as a generalization of a Sasakian structure in the direction where shear is present.
However, the shear is not arbitrary and is subject to relatively mild constraints.

\section{The $\eta$-Einstein condition on three-dimensional ACM manifolds}\label{sec:eta-einstein problem}
In this section, we discuss the extent to which three-dimensional normal ACM manifolds can admit an $\eta$-Einstein metric.
Here $\eta$-Einstein means that $Ric=ag+b\eta\otimes\eta,\ a,b\in C^\infty(M)$ holds.
The question of classifying three-dimensional normal ACM manifolds that are $\eta$-Einstein has long been a subject of interest.
In the non-compact case, the classification problem is difficult.
Even in the compact case, a classification result seems to be unavailable in the literature.
On the other hand, it has been shown that a three-dimensional contact metric manifold is $\eta$-Einstein if and only if it is a $(k,0,0)$-manifold \cite{Blair1990}.

When $\kappa=0$, the components of the Ricci tensor of an ACM structure are given by the generalized Sachs equations as
\begin{align}
&\textrm{Ric}(\xi,\xi)=-2\xi(\Theta)-2|\sigma|^2-2\Theta^2+2\omega^2,\\
&\textrm{Ric}(\partial,\partial)=P(\sigma)+2\Theta\sigma,\\
&\textrm{Ric}(\xi,\partial)=-\partial(\rho)-\bar{\eth}(\sigma).
\end{align}

From this we immediately obtain the following proposition.
\begin{screen}
\begin{prop}
A necessary and sufficient condition for a three-dimensional ACM structure with $\kappa = 0$ to be $\eta$-Einstein is that the following equations hold:
\begin{align}
P(\sigma) + 2\Theta\sigma = 0,\\
\partial(\rho) + \bar{\eth}(\sigma) = 0.
\end{align}
Both sides of the equations have the same spin weight, and hence the conditions are well-defined independently of the choice of the local adapted frame.
\end{prop}
\end{screen}

By applying the above proposition to three-dimensional contact metric structures, we recover the classification result for $\eta$-Einstein structures due to Blair et al. \cite{Blair1990}.
\begin{screen}
\begin{prop}
If a three-dimensional contact metric structure is $\eta$-Einstein, then it is a $(k,0,0)$-structure.
In this case $\textrm{Ric}(\xi,\xi)=2k$ and $|\sigma|$ are constants.
\end{prop}
\end{screen}
\begin{proof}
Since $\rho=i$, we have $P(\sigma)=\bar{\eth}(\sigma)=0$.
In this case we also have $2k:=\textrm{Ric}(\xi,\xi)=2(1-|\sigma|^2)$ from the above generalized Sachs equation.
By Proposition \ref{prop:(k,mu,nu)-mfd condition} it follows that the structure is a $(k,0,0)$-structure.
From the second Bianchi identity \eqref{eq:2nd bianchi 1}, we have $\partial k=0$.
Applying $[\partial,\bar{\partial}] = 2i\xi + \bar{\beta}\partial - \beta\bar{\partial}$ to $k$, it follows that $\xi(k)=0$ thus $k$ is a constant.
Therefore $|\sigma|^2$ is a constant.
\end{proof}

The condition corresponding to the $\eta$-Einstein condition $\partial\rho=0$ in three-dimensional normal ACM structures has been reported by Olszak \cite{Olszak}.
Compared with contact metric structures, normal ACM structures are shear-free, but the restrictions on the expansion $\Theta$ and the twist $\omega$ are weaker.
The function $\bar{\rho}$ is a CR function, but since we have $\rho=\Theta+i\omega$ and $\xi(\omega)+2\Theta\omega=0$, the real and imaginary parts are not independent functions.
Therefore one cannot directly apply the known results on three-dimensional CR functions.
In the $\eta$-Einstein problem for contact metric structures the main issue was to control the shear, whereas in the case of normal ACM structures the main issue is to control the expansion and the twist.

We approach the $\eta$-Einstein condition for three-dimensional compact normal ACM structures by deriving a priori estimates for the equation $\partial(\rho)+\bar{\eth}(\sigma)=0$ on a compact manifold.
\begin{screen}
\begin{prop}
On a three-dimensional compact ACM structure satisfying $\kappa=0$, if $\partial(\rho)+\bar{\eth}(\sigma)=0$ holds, then we have the following estimates:
\begin{align}
&\int_M\left(\frac{1}{2}||\nabla^h\Theta||^2+4\omega^2\Theta^2\right)dv\le
\int_M|\bar{\eth}\sigma|^2dv,\\
&\int_M\left(\frac{1}{2}||\nabla^h\omega||^2+4\omega^2\Theta^2\right)dv\le
\int_M|\bar{\eth}\sigma|^2dv,
\end{align}
where $\nabla^h$ denotes the horizontal gradient.
\label{prop:main inequality}
\end{prop}
\end{screen}

\begin{proof}
Let $U$ be a neighborhood and $\{\xi,\partial,\bar{\partial}\}$ be a local adapted complex frame on $U$.
We have
\begin{align}
&\partial \rho = -(\overline{\partial}\sigma + 2\sigma\overline{\beta}),\\
&\frac{1}{\sqrt{2}}(e_2 - i e_3)(\Theta + i\omega) = -(\overline{\partial}\sigma + 2\sigma\overline{\beta})=:A+iB,\\
&\frac{1}{\sqrt{2}} \{ (e_2\Theta + e_3\omega) + i(e_2\omega - e_3\Theta) \} = A + iB,\\
\therefore\ &e_2\Theta = -e_3\omega + \sqrt{2}A,\\
&e_3\Theta = e_2\omega - \sqrt{2}B.
\end{align}
which implies
\begin{align}
\nabla^h \Theta &= (-e_3\omega + \sqrt{2}A)e_2 + (e_2\omega - \sqrt{2}B)e_3\\
&= (-e_3\omega)e_2 + (e_2\omega)e_3 + \sqrt{2}(Ae_2 - Be_3)\\
&=\varphi\nabla^h\omega+F_\sigma,\\
F_\sigma &= -\sqrt{2} \left( \text{Re}(\overline{\partial}\sigma + 2\sigma\overline{\beta}) e_2 - \text{Im}(\overline{\partial}\sigma + 2\sigma\overline{\beta}) e_3 \right),\\
||F_\sigma||^2 &= 2 |\overline{\partial}\sigma + 2\sigma\overline{\beta}|^2=2|\bar\eth\sigma|^2.
\end{align}

Next we calculate $\textrm{div}(\varphi \nabla^h \omega)$ so that we derive a subelliptic operator $\Delta^h\Theta=\textrm{div}(\nabla^h\Theta)$.
Since we have $\varphi \nabla^h \omega=-e_3(\omega){e_2}+{e_2}(\omega)e_3$, the following holds.
\begin{align}
\text{div}(\varphi \nabla^h \omega) &= \text{div}\bigl( ({e_2}\omega)e_3 \bigr) - \text{div}\bigl( (e_3\omega){e_2} \bigr)\\
&= e_3({e_2}\omega) + ({e_2}\omega)\text{div}e_3-{e_2}(e_3\omega) - (e_3\omega)\text{div}{e_2}\\
&= -[{e_2},e_3]\omega  + ({e_2}\omega)\text{div}e_3 - (e_3\omega)\text{div}{e_2}.
\end{align}

Here we have
\begin{align}
i[{e_2},e_3]&=[\partial,\bar{\partial}]=2i\omega\xi+\bar{\beta}\partial-\beta\bar{\partial},\\
[{e_2},e_3]&=2\omega\xi-i(\bar{\beta}\partial-\beta\bar{\partial}),
\end{align}
and
\begin{align}
&g([{e_2},e_3],{e_2})=g(\nabla_{e_2}e_3,{e_2})-g(\nabla_{e_3}{e_2},{e_2}) =g(\nabla_{e_2}e_3,{e_2})=\textrm{div}(e_3),\\
&g([{e_2},e_3],e_3)=g(\nabla_{e_2}e_3,e_3)-g(\nabla_{e_3}{e_2},e_3)=-g(\nabla_{e_3}{e_2},e_3)=-\textrm{div}({e_2}),
\end{align}
which implies
\begin{align}
[{e_2}, e_3] = 2\omega\xi + (\text{div}e_3){e_2} - (\text{div}{e_2})e_3.
\end{align}
Therefore we obtain
\begin{align}
\text{div}(\varphi \nabla^h \omega) &=
-2\omega\xi\omega - (\text{div}e_3){e_2}\omega + (\text{div}{e_2})e_3\omega
+ ({e_2}\omega)\text{div}e_3 - (e_3\omega)\text{div}{e_2}
=-2\omega\xi\omega\\
&=4\Theta\omega^2\ (\because \xi\omega = -2\Theta\omega).
\end{align}
It follows that
\begin{align}
\text{div}(\nabla^h \Theta)
&=\textrm{div}(\varphi \nabla^h \omega) + \textrm{div}(F_\sigma),\\
\Delta^h \Theta-4\omega^2& \Theta=\textrm{div}(F_\sigma).
\end{align}
This equation holds on the neighborhood $U$, but since $\rho=\Theta+i\omega$ is independent of the choice of the local adapted frame, so is the left-hand side.
Moreover, one can verify by a direct calculation that $F_\sigma = -\sqrt{2} \left( \operatorname{Re}(\bar\eth\sigma) e_2 - \operatorname{Im}(\bar\eth\sigma) e_3 \right)$ is also independent of the choice of the local adapted frame.
Indeed, when $\partial'=e^{i\theta}\partial$ and $(\bar\eth\sigma)'=e^{i\theta}\bar\eth\sigma$, we have
\begin{align}
 \operatorname{Re}(\bar\eth\sigma) e_2 - \operatorname{Im}(\bar\eth\sigma) e_3
 &=(e_2,e_3)
 \begin{pmatrix}
\operatorname{Re}(\bar\eth\sigma)\\
- \operatorname{Im}(\bar\eth\sigma)
 \end{pmatrix}\\
&=(e_2,e_3)
\begin{pmatrix}
\cos\theta & -\sin\theta\\
\sin\theta & \cos\theta
\end{pmatrix}
\begin{pmatrix}
\cos\theta & \sin\theta\\
-\sin\theta & \cos\theta
\end{pmatrix}
 \begin{pmatrix}
\operatorname{Re}(\bar\eth\sigma)\\
- \operatorname{Im}(\bar\eth\sigma)
 \end{pmatrix}\\
&=(e'_2,e'_3)
 \begin{pmatrix}
\operatorname{Re}(\bar\eth'\sigma')\\
- \operatorname{Im}(\bar\eth'\sigma')
 \end{pmatrix}
=\operatorname{Re}(\bar\eth'\sigma') e'_2
-\operatorname{Im}(\bar\eth'\sigma') e'_3 .
\end{align}
Therefore, the right-hand side is also independent of the choice of the local adapted frame, and hence this equation holds on all of $M$.

Thus $\Theta$ is a solution of a nonlinear subelliptic equation whose source term is $\textrm{div}(F_\sigma)$, where there is a nonlinear dependence between $\Theta$ and $\omega$.
Multiplying this equation by $\Theta$ and integrating, we obtain
\begin{align}
\int_M(||\nabla^h\Theta||^2+4\omega^2\Theta^2)dv=-\int_M\Theta \textrm{div}(F_\sigma)dv=\int_Mg(\nabla^h\Theta,F_\sigma)dv .
\end{align}
Using Young's inequality $ab \le \epsilon a^2 + \frac{1}{4\epsilon}b^2,\ \epsilon=1/2$, we obtain the estimate
\begin{align}
\int_M\left(\frac{1}{2}||\nabla^h\Theta||^2+4\omega^2\Theta^2\right)dv\le\frac{1}{2}\int_M||F_\sigma||^2dv= \int_M|\bar{\eth}\sigma|^2dv.
\end{align}
Moreover, using $\nabla^h \Theta = \varphi \nabla^h \omega + F_\sigma$, rewriting the above in terms of $\nabla^h\omega$, and applying the same type of estimate, we obtain
\begin{align}
\frac{1}{2} \int_M ||\nabla^h \omega||^2 dv + \int_M 4\Theta^2 \omega^2 dv \le  \int_M |\overline{\eth}\sigma|^2 dv .
\end{align}
\end{proof}

From the above estimates we immediately obtain the following result.
\begin{screen}
\begin{thm}
If a compact three-dimensional normal ACM manifold, or equivalently trans-Sasakian manifold is $\eta$-Einstein, then the structure is either an $\alpha$-Sasakian structure homothetic to a Sasakian structure or a $\beta_s$-Kenmotsu structure of expansion $\Theta$ satisfying $\partial\Theta=0$.
\end{thm}
\end{screen}
\begin{proof}
Since $\kappa=\sigma=0$, Proposition \ref{prop:main inequality} and \eqref{eq:1st bianchi 1} give $\partial\omega=\Theta\omega=0$ and $\xi(\omega)=-2\Theta\omega=0$.
Thus $\omega$ is constant on each connected component, so the rank is constant.
In the case of rank one, we have $\omega=0$ and $\partial\Theta=0$ from the inequality of Proposition \ref{prop:main inequality}, which implies that the structure is a $\beta$-Kenmotsu structure.
In the case of rank three, we have $\Theta=0$ and $\partial\omega=0$.
Since $[\partial,\bar{\partial}] = 2i\omega\xi + \bar{\beta}\partial - \beta\bar{\partial}$, it follows that $\xi(\omega)=0$ thus $\omega$ is constant.
If necessary, by reversing the orientation of $\xi$, we may assume that $\omega=c>0$.
Consider the horizontal conformal transformation $g'=c\,g+(1-c)\eta\otimes\eta$.
Then $g'$ is a Sasakian metric.
Applying the $D$-homothetic transformation with factor $c$ to the Sasakian metric $g'$, we obtain another Sasakian metric $g''=c\,g'+(c^2-c)\eta\otimes\eta$.
Thus $g''=c^2g$ holds.
\end{proof}

\begin{rem}
There exists an example, known as the Rossi sphere, which is of rank $3$, compact, and $\eta$-Einstein with respect to the Webster curvature, but is not $\alpha$-Sasakian. 
However, this example is not $\eta$-Einstein with respect to the usual Ricci tensor, and since $\sigma\neq 0$ it is not normal \cite{Ho}.
\end{rem}

\section{Conclusions and Discussion}
\begin{figure}[H]
\includegraphics[scale=0.4]{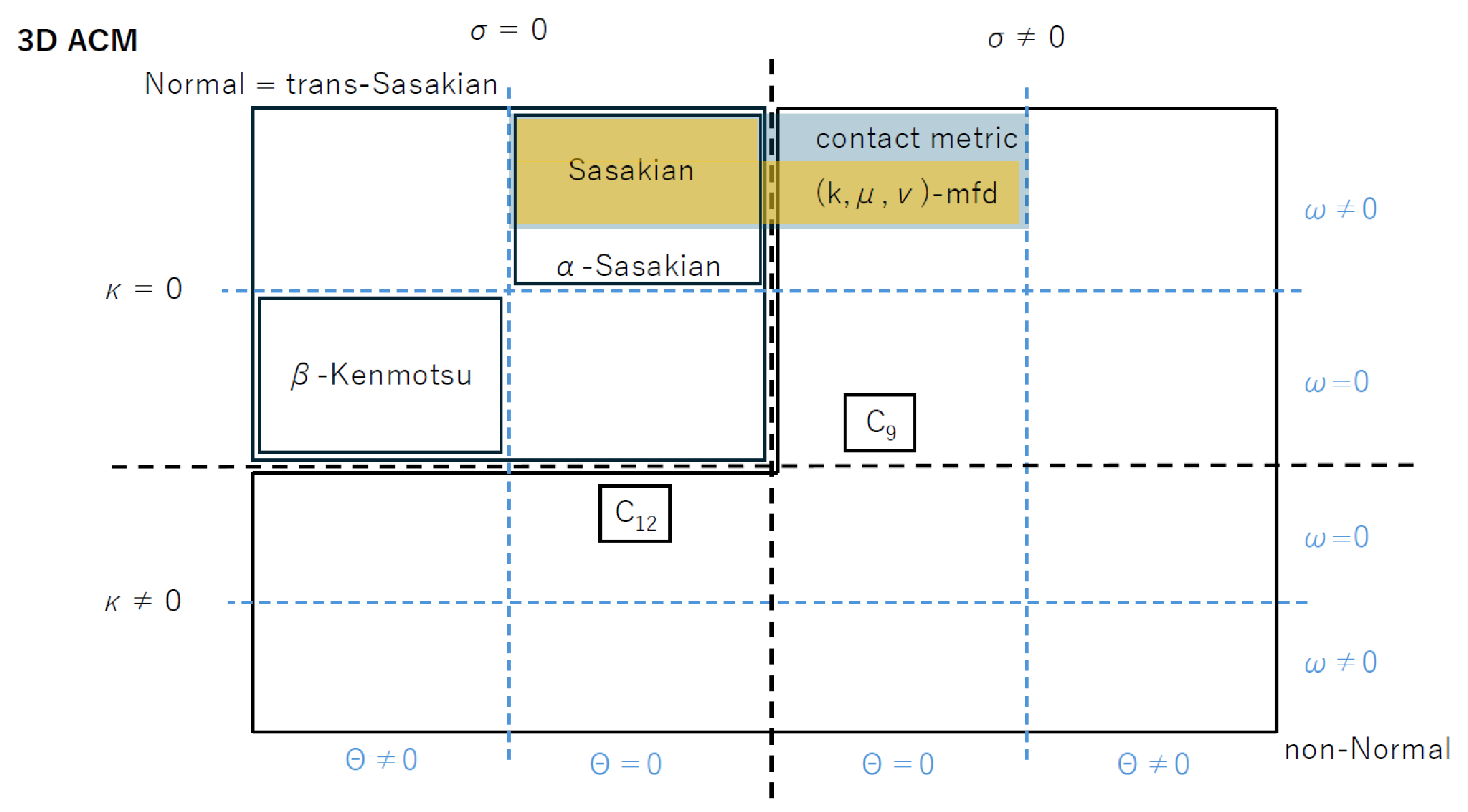}
\caption{Relationship between traditional classes and the partition of three-dimensional ACM structures based on $\Theta, \omega$, and $\sigma$}
\label{fig:3D ACM map}
\end{figure}

In this paper, we have discussed three-dimensional ACM structures in terms of spin coefficients. The spin coefficients particularly relevant to the characterization of traditional classes of ACM structures are $\kappa$, the expansion $\Theta$, the twist $\omega$, and the complex shear $\sigma$.
Figure \ref{fig:3D ACM map} illustrates the relationship between the classes discussed in this paper—namely, trans-Sasakian structures, contact metric structures, and $(k,\mu,\nu)$-structures—and the partition of three-dimensional ACM structures based on $\kappa, \sigma, \Theta,$ and $\omega$.
Furthermore, the well-known classes of three-dimensional non-normal ACM manifolds, $C_9$ and $C_{12}$, are also positioned within the non-normal region.
Traditional definitions of ACM subclasses can include algebraic normalization conditions that are of little relevance to their analytic properties.
The approach using spin coefficients allows for a flexible treatment by relaxing such conditions.
Thus, the Newman-Penrose formalism demonstrates its utility for analytic properties.

% Finally, we comment on the classification of three-dimensional normal ACM $\eta$-Einstein manifolds.
% Although we have provided a classification for the compact case in this study, the non-compact case has been left untouched.
% Since, in the case of rank 1, the construction of $\beta$-Kenmotsu $\eta$-Einstein structures in the non-compact setting is analogous to the compact one, the non-trivial case of interest is rank 3.
% By using the Raychaudhuri equations obtained from the generalized Sachs equation,
% $$-\xi(\Theta) = \Theta^2 - \omega^2 + \frac{1}{2}\textrm{Ric}(\xi,\xi), \quad -\xi(\omega) = 2\Theta\omega,$$
% the behavior of $\omega$ and $\Theta$ along the Reeb orbits can be fully described in terms of elementary functions.
% However, a major difficulty arises in the non-compact case regarding the control of these quantities in the horizontal directions.
% In particular, when $\inf\omega=0$ or the $\varphi$-sectional curvature is not bounded below, $\eta$-Einstein solutions are expected to exhibit complex behavior even in the complete case.

\section*{Acknowledgement}
The author would like to thank Ueno for pointing out several computational errors.
The author is also grateful to Ishihara, Koike, Kozaki, and Morisawa for their helpful comments and suggestions.
This work was partly supported by MEXT Promotion of Distinctive Joint Research Center Program JPMXP0723833165 and Osaka Metropolitan University Strategic Research Promotion Project (Development of International Research Hubs).

\appendix
\section{Spin Covariance of $P$, $\eth$, and $\bar{\eth}$}\label{App: P eth spin covariance}
\begin{align*}
(P q)' &= \xi(q') - s \epsilon' q' \\
&= \xi(e^{is\theta}q) - s(\epsilon + i\xi(\theta)) (e^{is\theta}q) \\
&= \left[ i s \xi(\theta) e^{is\theta} q + e^{is\theta} \xi(q) \right] - \left[ s\epsilon e^{is\theta}q + i s \xi(\theta) e^{is\theta} q \right] \\
&= e^{is\theta} \xi(q) - s\epsilon e^{is\theta}q \\
&= e^{is\theta} (\xi(q) - s \epsilon q) \\
&= e^{is\theta} (P q)\\
(\eth q)' &= \partial'(q') - s \beta' q' \\
&= e^{i\theta}\partial (e^{is\theta}q) - s \left[ e^{i\theta}(\beta + i \partial \theta) \right] (e^{is\theta}q) \\
&= e^{i\theta} \left[ i s (\partial \theta) e^{is\theta} q + e^{is\theta} \partial(q) \right] - s e^{i(s+1)\theta}(\beta + i \partial \theta) q \\
&= i s (\partial \theta) e^{i(s+1)\theta} q + e^{i(s+1)\theta} \partial(q) - s \beta e^{i(s+1)\theta} q - i s (\partial \theta) e^{i(s+1)\theta} q\\
&= e^{i(s+1)\theta} \partial(q) - s \beta e^{i(s+1)\theta} q \\
&= e^{i(s+1)\theta} (\partial q - s \beta q) \\
&= e^{i(s+1)\theta} (\eth q)\\
(\overline{\eth} q)' &= \overline{\partial}'(q') + s \overline{\beta}' q' \\
&= e^{-i\theta}\overline{\partial} (e^{is\theta}q) + s \left[ e^{-i\theta}(\overline{\beta} - i \overline{\partial} \theta) \right] (e^{is\theta}q) \\
&= e^{-i\theta} \left[ i s (\overline{\partial} \theta) e^{is\theta} q + e^{is\theta} \overline{\partial}(q) \right] + s e^{i(s-1)\theta}(\overline{\beta} - i \overline{\partial} \theta) q \\
&= i s (\overline{\partial} \theta) e^{i(s-1)\theta} q + e^{i(s-1)\theta} \overline{\partial}(q) + s \overline{\beta} e^{i(s-1)\theta} q - i s (\overline{\partial} \theta) e^{i(s-1)\theta} q\\
&= e^{i(s-1)\theta} \overline{\partial}(q) + s \overline{\beta} e^{i(s-1)\theta} q \\
&= e^{i(s-1)\theta} (\overline{\partial} q + s \overline{\beta} q) \\
&= e^{i(s-1)\theta} (\overline{\eth} q)
\end{align*}

\end{document}